\documentclass[12pt]{article}
\usepackage[cp1251]{inputenc}
\usepackage[russian]{babel}
\usepackage{amssymb,amsthm,amsmath}

\theoremstyle{plain}

\begin{document}

{\large {\bf \centerline{Асимптотические формулы для собственных значений}}
{\bf \centerline{ и собственных функций
операторов Штурма--Лиувилля}} {\bf \centerline{ с потенциалами ---
распределениями.}}
{\bf \centerline{Краевые условия Дирихле--Неймана. }}

\bigskip

{\centerline{Швейкина~О.~А.}}

\medskip
В настоящей статье изучается оператор Штурма--Лиувилля
\begin{equation}\label{1}
Ly=-\dfrac{d^2y}{dx^2}+q(x)y,
\end{equation}
в пространстве $L_2[0,\pi]$ с граничными условиями
Дирихле--Неймана $y(0)=y'(\pi)=0$. Предполагается, что потенциал
$q(x)=u'(x),$ где $u\in L_2[0,\pi]$. Производная здесь понимается
в смысле распределений. Операторы такого (и более общего) вида
были определены в работах \cite{SSh1}--\cite{SSh2}. Там же было
доказано, что оператор $L$ имеет чисто дискретный спектр, а в
случае вещественного потенциала
--- самосопряжен и полуограничен. При этом в формулах для
собственных значений были выписаны два первых члена асимптотики;
для собственных функций был получен главный член. В недавней
статье А.М.Савчука \cite{S1} была доказана более подробная
асимптотика для собственных функций в случае краевых условий
Дирихле. Полученные асимптотические формулы имеют приложения к
задачам разного характера, например, в работах И.В.Садовничей
\cite{Sad1}--\cite{Sad2}, эти результаты использовались при
доказательстве теорем равносходимости.

Цель данной статьи --- получить асимптотические формулы для
собственных значений и собственных и присоединенных функций
оператора \eqref{1} с краевыми условиями Дирихле--Неймана.
Потенциал $q(x)$ предполагается комплекснозначным. Некоторые
результаты об операторе \(L\), полученные в \cite{SSh2} и
необходимые в этой работе, будут приведены ниже.

Наша задача --- получить асимптотические формулы для решений
уравнения
\begin{equation} \label{2b} -y'' + q(x)y = \lambda
y.\end{equation}

 Уравнение \eqref{2b} можно записать в виде
системы (см. \cite{SSh2})
\begin{equation}\label{3b}\left( \begin{matrix} y_1 \\ y_2
\end{matrix}\right)' = \left( \begin{matrix} u & 1 \\ -\lambda-u^2
& -u \end{matrix}\right)\left( \begin{matrix} y_1 \\ y_2
\end{matrix}\right), \ \ y_1 = y,\  y_2 = y' - u(x)y.\end{equation}
Сделаем замену \(y_1(x, \lambda) = r(x, \lambda)\sin{\theta(x,
\lambda)}, \ y_2(x, \lambda) = \lambda^{\frac{1}{2}}r(x,
\lambda)\cos{\theta(x, \lambda)},\) которая является модификацией
замены Прюфера (см. \cite{H}). Тогда систему \eqref{3b} можно
переписать в виде
\begin{equation} \label{2} \begin{aligned} r'\sin{\theta} +
r\theta'\cos{\theta} &=
ur\sin{\theta}+\lambda^{\frac{1}{2}}r\cos{\theta} \\
\lambda^{\frac{1}{2}}r'\cos{\theta}-\lambda^{\frac{1}{2}}r\theta'\sin{\theta}
&= -\lambda r \sin{\theta} - u^2 r \sin{\theta} -
\lambda^{\frac{1}{2}}ur\cos{\theta},\\ \end{aligned}\end{equation}
где \(r = r(x, \lambda),\ \theta = \theta(x, \lambda), u = u(x)\),
а производные функций \(r\) и \(\theta\) берутся по переменной
\(x\). Умножим первое уравнение в \eqref{2} на
\(\lambda^{\frac{1}{2}}\cos{\theta}\) и вычтем второе уравнение,
умноженное на \(\sin{\theta}\). В результате получим уравнение для
функции \(\theta(x, \lambda)\)
 \begin{equation}
 \label{3}
 \theta'(x,
 \lambda) = \lambda^{1/ 2} + \lambda^{-1/ 2}u^2(x)\sin^2{\theta(x, \lambda)}+u(x)\sin{2\theta(x, \lambda)}.
 \end{equation}
Если мы сложим первое уравнение в \eqref{2}, умноженное на
\(\lambda^{\frac{1}{2}}\sin{\theta}\), со вторым уравнением,
умноженным на \(\cos{\theta}\), то получим уравнение на функцию
\(r(x, \lambda)\) \begin{equation} \label{4} r'(x, \lambda) =
-r(x, \lambda)\left[u(x)\cos{2\theta(x,
\lambda)}+\frac{1}{2}\lambda^{-1/ 2}u^2(x)\sin{2\theta(x,
\lambda)}\right].
\end{equation}

Итак, мы получили формулы для функций $r(x,\lambda)$ и
$\theta(x,\lambda)$, с помощью которых можно выразить решение
уравнения \eqref{2b}. Введем необходимые обозначения и
сформулируем вспомогательные леммы для исследования этих формул.
Будем обозначать
\begin{equation}\begin{aligned}\label{7a}
v(x,\lambda) = &\int_0^x u(t)\sin(2\lambda^{1/ 2}t)dt +
\frac{1}{2}\lambda^{-1/2}\int_0^x u^2(t)dt+ \\
&+2\int_0^x\int_0^t u(t)u(s)\cos(2\lambda^{1/
2}t)\sin(2\lambda^{1/ 2}s)dsdt- \\ &-\frac{1}{2}\lambda^{-1/
2}\int_0^x u^2(t)\cos(2\lambda^{1/ 2}t)dt,\\
\end{aligned}
\end{equation}

\[\begin{aligned}
\gamma(\lambda) = & \sup_{0\leq x \leq \pi}\biggl(\left|\int_0^x
u(t)\sin (2\lambda^{1/ 2}t) dt\right| + \left|\int_0^x u(t)\cos
(2\lambda^{1/ 2}t) dt \right|+ \\ & + 2\left|\int_0^x\int_0^t
u(t)u(s)\cos(2\lambda^{1/ 2}t) \sin(2\lambda^{1/ 2}s)dtds\right| +
\\ &
 +\frac{1}{2}\left|\lambda^{-1/ 2}\int_0^x u^2(t)\cos (2\lambda^{1/ 2}t)  dt \right|\biggr) +
  |\lambda|^{-1/ 2}||u||_{L_2}^2.\\
\end{aligned}\]

Сформулируем утверждения, которые потребуются нам в дальнейшем.

\newpage

 \textbf{Лемма 1} (см. \cite{SSh2})

 Пусть \(\alpha > 0\) - произвольное фиксированное число, а \(P_\alpha\) - область, ограниченная параболой
 \(|\text{Im}\  \sqrt{\lambda}| < \alpha\). Тогда существует число
 \(\mu\), зависящее только от \(u(x)\) и \(\alpha\) такое, что при
 любых  \(\lambda \in P_\alpha\),
 \(\text{Re}\ {\lambda} > \mu \) уравнение
\eqref {3} имеет единственное решение \(\theta(x, \lambda)\),
определенное при всех \(0\leq x \leq \pi\) и удовлетворяющее
начальному условию \(\theta(0, \lambda) = 0\). Это решение
допускает представление
\[\theta(x, \lambda) = \lambda^{1/ 2}x + v(x,
\lambda) + \rho(x, \lambda),\] где \( \sup\limits_{0\leq x \leq
\pi}{|\rho(x, \lambda)|}\leq M\gamma^2(\lambda),\  \lambda \in
P_\alpha,\ \text{Re}\ {\lambda} > \mu,\) причем выбор постоянной
\(M\) зависит от функции \(u\) и \(\alpha\), но не зависит от
\(x\), \(\lambda\).

\textbf{Лемма 2} (см. \cite{SSh2})

Пусть \(\alpha > 0\) - произвольное фиксированное число, а
\(P_\alpha\) - область, ограниченная параболой
 \(|\text{Im}\  {\sqrt{\lambda}}| < \alpha\).
 Пусть \(\theta(x, \lambda)\) - решение уравнения \eqref{3}
  с начальным
условием \(\theta(0, \lambda)\) = 0. Тогда решение $r(x, \lambda)$
уравнения \eqref{4} с начальным условием $r(0, \lambda) = 1$
допускает представление
\[r(x, \lambda) = 1 - \int_0^x{u(t)\cos{(2nt)}dt} -
\frac{1}{2n}\int_0^x{u^2(t)\sin{(2nt)}dt} + \rho(x, \lambda),\]\[
\lambda \in P_\alpha,\ \text{Re}\  {\lambda} > \mu,\] где
\(\sup\limits_{0\leq x \leq \pi}{|\rho(x, \lambda)|}\leq
M\gamma^2(\lambda)\), причем \(M\) и \(\mu\) зависят только от
\(u\) и \(\alpha\).

Переходим к формулировке и доказательству утверждений
непосредственно данной статьи.

\textbf{Теорема 1.} Для собственных значений оператора \(L = -
d^2/ dx^2 - q(x)\) с граничными условиями Дирихле-Неймана \(y(0) =
0, \ y'(\pi) = 0\), где \(q(x) = u'(x)\), а \( u(x) \in L_2\),
выполнено:

\begin{equation}\label{7b}\lambda_n^{1/ 2} = n-\frac{1}{2} -\frac{1}{\pi}v(\pi,
(n-\tfrac{1}{2})^2) +\rho(\lambda_n),\end{equation} где через
\(\rho(\lambda_n)\) здесь и далее будем обозначать произвольную
последовательность, удовлетворяющую условию \(
|\rho(\lambda_n)|\leq M\gamma^2(\lambda_n).\) Здесь \(M\) зависит
только от функции $u$.

\textbf{Доказательство}

Обозначим \(m = n-\frac{1}{2}, n \in \mathbb{N}\). Для уравнения
\(-y'' = \lambda y\) общее решение имеет вид \(y =
c_1\sin{(\lambda^{1/2}x)}+c_2\cos{(\lambda^{1/2}x)}\), где \(c_1\)
и \(c_2\) --- произвольные постоянные. Подставив в решение
граничные условия, получим, что

1. \(y(0) = 0\), откуда следует, что \(c_2 = 0\);

2. \(y'(\pi) = 0\), значит, \(\cos{\theta(\pi, \lambda_n)} = 0\),
и следовательно, \(\theta(\pi, \lambda_n) = \pi m\). По Лемме 1
при \(x = \pi\) получим:
\[\lambda_n^{1/ 2} + \frac{1}{\pi}v(\pi, \lambda_n) + \rho(\lambda_n) = m,
\ \text{где}\  |\rho(\lambda_n)|\leq M\gamma^2(\lambda_n).\]
 Таким образом, необходимо доказать оценку
\begin{equation}\label{6}|v(\pi, \lambda_n) - v(\pi, m^2)| \leq
M\gamma^2(\lambda_n).\end{equation}
Преобразуем левую часть
неравенства.
\begin{equation}\label{7}\begin{aligned}
v(\pi, \lambda_n) - v(\pi, m^2) =
\int_0^{\pi}{u(t)(\sin{(2\lambda_n^{1/ 2}t)} - \sin{(2mt)})dt} +
\\ + \frac{1}{2} (\lambda_n^{-1/ 2}-m^{-1})\left(\int_0^{\pi}{u^2 (t)dt} -
\int_0^{\pi}{u^2(t)\cos{(2\lambda_n^{1/ 2}t)}dt} \right)-\\ -
\frac{1}{2} m^{-1}\left(\int_0^{\pi}{u^2(t)\cos{(2\lambda_n^{1/
2}t)}dt} - \int_0^{\pi}{u^2(t)\cos{(2mt)}dt} \right)+
\\ +2\left(\int_0^\pi\int_0^t u(t)u(s)\cos (2\lambda_n^{1/
2}t)\sin(2\lambda_n^{1/ 2}s)dsdt- \right.
\\- \left.\int_0^\pi\int_0^t u(t)u(s)\cos(2mt)\sin(2ms) ds dt
\right)\\
\end{aligned}
\end{equation}

Обозначим слагаемые из правой части \eqref {7} в порядке их
очередности \(I_1, I_2, I_3, I_4\), и введем величину \(\nu_n =
\lambda_n^{1/ 2} - m\), причем отметим, что из \eqref {7a}
следует, что \(|\nu_n|\leq M \gamma^2(\lambda_n)\). Перейдем к
оценке \(I_1\), для этого преобразуем сумму синусов к более
удобному виду.

\[\begin{aligned}
&\sin{(2\lambda_n^{1/ 2}t)} - \sin{(2mt)} = \sin{(2\lambda_n^{1/
2}t)} - \sin{(2\lambda_n^{1/ 2}t - 2\nu_n t)}
=\\
=&\sin{(2\lambda_n^{1/ 2}t)} - \sin{(2\lambda_n^{1/
2}t)}\cos{(2\nu_n t)} + \cos{(2\lambda_n^{1/ 2}t)}\sin{(2\nu_n t)}
=\\ =&\sin{(2\lambda_n^{1/ 2}t)} - \sin{(2\lambda_n^{1/ 2}t)}(1 -
O(\nu_n ^ 2)) + \cos{(2\lambda_n^{1/ 2}t)}(2\nu_n t + O(\nu_n ^
3)) = \\ =&2\nu_n t\cos{(2\lambda_n^{1/ 2}t)} + O(\nu_n ^ 2)\\
\end{aligned}\]
 Подставим полученное в \(I_1\):
\[ \begin{aligned} I_1 &=
\left|2\nu_n\int_0^{\pi}{tu(t)\cos{(2\lambda_n^{1/ 2}t)}dt} +
O(\nu_n ^ 2)\right|= \bigg
|2\nu_n\pi\int_0^{\pi}{u(t)\cos{(2\lambda_n^{1/ 2}t)}dt}
\\  &- 2\nu_n\int_0^{\pi}\int_0^t{u(s)\cos{(2\lambda_n^{1/
2}s)}dsdt}\bigg| + O(\nu_n ^ 2)\leq M |\nu_n| |\gamma(\lambda_n)|
\leq M^2 \gamma^2(\lambda_n).\\ \end{aligned}\]
 Далее заметим, что
 \(\lambda_n^{-1/ 2}-m^{-1} = O(m^{-2})\), и \(\cos{(2\lambda_n^{1/
2}t)} - \cos{(2mt)} = O(\nu_n)\), причем \(m^{-2}\leq M
\gamma^2(\lambda_n) \) и \(\nu_n m^{-1}\leq M
\gamma^2(\lambda_n)\), следовательно \(I_2\) и  \(I_3\) по модулю
не превосходят величины \(M \gamma^2(\lambda_n)\).

Для оценки последнего слагаемого в правой части \eqref{7}
воспользуемся равенством:
 \[\begin{aligned}&\cos{(2\lambda_n^{1/
2}t)}\sin{(2\lambda_n^{1/ 2}s)} - \cos{(2mt)}\sin{(2ms)} = \\ =
&\cos{(2\lambda_n^{1/ 2}t)}(\sin{(2\lambda_n^{1/
2}s)}-\sin{(2ms)}) + (\cos{(2\lambda_n^{1/
2}t)}-\cos{(2mt)})\sin{(2ms)} = \\ = &2\nu_n s
\cos{(2\lambda_n^{1/ 2}t)}\cos{(2\lambda_n^{1/ 2}s)} -
\sin{(2\lambda_n^{1/ 2}t)}\sin{(2\lambda_n^{1/ 2}s)} + O(\nu_n ^
2)\end{aligned}\]

Из определения функции $\gamma(\lambda_n)$ после изменения порядка
интегрирования получим:
\[2\nu_n\left(\int_0^{\pi}{u(t)\cos{(2\lambda_n^{1/
2}t)}\int_0^t{u(s)\cos{(2\lambda_n^{1/ 2}s)}dsdt}}\right)\leq M
\gamma^2(\lambda_n).\] То же верно для аналогичного интеграла с
синусами. Тогда \(|I_4| \leq M \gamma^2(\lambda_n)\), и
справедливость оценки \eqref{6} установлена.

\emph{Утверждение доказано.}

\textbf{Замечание}

Из \cite{SSh2} (Лемма 2.4) следует, что в случае \(u(x) \in L_2\)
- последовательность \( \{ \gamma(\lambda_n)\} \in l_2\), и таким
образом, \( \{ \gamma^2(\lambda_n)\} \in l_1\).

Переходим к вопросу о собственных и присоединенных функциях
оператора \eqref{1}.

\textbf{Теорема 2}

Рассмотрим оператор \(L\), порожденный дифференциальным выражением
\(-y'' + q(x)y \), где \(q(x) = u'(x)\) в смысле распределений, а
комплекснозначная функция \(u(x) \in L_2\), и краевыми условиями
Дирихле-Неймана \(y(0) = 0, \ y'(\pi) = 0\). Обозначим через
\(\{y_n(x)\}_{n = 1}^\infty\) систему собственных и присоединенных
функций оператора \(L\) (определение системы собственных и
присоединенных функций см. в \cite{S1}), через \(\{w_n(x)\}_{n =
1}^\infty\)
--- биортогональную систему. Тогда справедлива следующая
асимптотическая формула:

\begin{multline} \label{8} \sqrt{\frac{\pi}{2}}\ y_n(x) =\\
\sin(mx)\left(1 + \frac{1}{\pi}\int_0^\pi (\pi -
t)u_R(t)\cos(2mt)dt - \int_0^x u(t)\cos(2mt)dt\right) +
\\ +\frac{1}{2m}\sin(mx)\biggl( -\int_0^x u^2(t)\sin(2mt)dt
+\frac{1}{\pi}\int_0^\pi (\pi - t)(u_R^2(t) - u_I^2(t)\sin(2mt)dt)
\biggr) + \\  + \cos(mx)\biggl(\int_0^x
u(t)\sin(2mt)dt+2\int_0^x\int_0^t u(t)u(s)\cos(2mt)\sin(2ms)dsdt-
\\  -\frac{x}{\pi}\int_0^\pi
u(t)\sin(2mt)dt-\frac{2x}{\pi}\int_0^\pi\int_0^t
u(t)u(s)\cos(2mt)\sin(2ms)dsdt \biggr) +\\  +
\frac{1}{2m}\cos(mx)\biggl(\int_0^x u^2(t)dt - \int_0^x
u^2(t)\cos(2mt)dt-\\  -\frac{x}{\pi}\int_0^\pi u^2(t)dt +
\frac{x}{\pi}\int_0^\pi u^2(t)dt\cos(2mt) \biggr) +\rho(x,
\lambda_n),
\end{multline}

Соответствующие функции биортогональной системы имеют вид:
\begin{multline} \label{9b} \sqrt{\frac{\pi}{2}}\ w_n(x) =\\
 \sin(mx)\left(1 + \frac{1}{\pi}\int_0^\pi (\pi -
t)(u_R(t)+2iu_I(t))\cos(2mt)dt - \int_0^x
\overline{u}(t)\cos(2mt)dt\right) +
\\ +\frac{1}{2m}\sin(mx)\biggl( -\int_0^x
\overline{u}^2(t)\sin(2mt)dt+\\ +\frac{1}{\pi}\int_0^\pi (\pi -
t)(u_R^2(t) - u_I^2(t)+4iu_R(t)u_I(t))\sin(2mt)dt) \biggr) + \\
+ \cos(mx)\biggl(\int_0^x
\overline{u}(t)\sin(2mt)dt+2\int_0^x\int_0^t
\overline{u}(t)\overline{u}(s)\cos(2mt)\sin(2ms)dsdt-
\\  -\frac{x}{\pi}\int_0^\pi
\overline{u}(t)\sin(2mt)dt-\frac{2x}{\pi}\int_0^\pi\int_0^t
\overline{u}(t)\overline{u}(s)\cos(2mt)\sin(2ms)dsdt \biggr) +\\
+ \frac{1}{2m}\cos(mx)\biggl(\int_0^x \overline{u}^2(t)dt -
\int_0^x \overline{u}^2(t)\cos(2mt)dt-\\
-\frac{x}{\pi}\int_0^\pi \overline{u}^2(t)dt +
\frac{x}{\pi}\int_0^\pi \overline{u}^2(t)dt\cos(2mt) \biggr)
+\rho(x, \lambda_n),
\end{multline}

 где \(\{ \sup\limits_{0 \leq x \leq \pi}{|\rho(x,
\lambda_n)|} \} \in l_1\), \(m = n - \frac{1}{2}, n \in
\mathbb{N}\), \(u_R(x)\) и \(u_I(x)\) обозначают вещественную и
мнимую части функции \(u(x)\) соответственно.

\textbf{Доказательство}

Для вычисления собственных и присоединенных функций будем
использовать формулу
\begin{equation}\label{10}y_n (x, \lambda_n) = r(x, \lambda_n)\sin{\theta(x,
\lambda_n)} ,\end{equation} справедливость которой в рамках данной
задачи показана в статье \cite{SSh2}, где была проделана
аналогичная работа с краевыми условиями Дирихле. Подставим в
\eqref{10} асимптотические формулы для функций \(r(x, \lambda)\) и
\(\theta(x,\lambda)\), где согласно Леммам 1 и 2

\begin{equation*}\begin{aligned}&r(x, \lambda_n) = 1 - \int_0^x{u(t)\cos{(2\lambda_n^{1/
2}t)}dt} - \frac{1}{2}\lambda_n^{-1/
2}\int_0^x{u^2(t)\cos{(2\lambda_n^{1/ 2}t)}dt} + \rho(x,
\lambda_n),\\ &\theta(x, \lambda_n) = \lambda_n^{1/ 2}x + v(x,
\lambda_n) + \rho(x, \lambda_n), \sup_{0\leq x \leq \pi}{|\rho(x,
\lambda_n)|}\leq M \gamma^2(\lambda_n)\end{aligned}\end{equation*}

Обозначим \(\mu_n = -\frac{1}{\pi}v(\pi, m^2).\)  Далее
воспользуемся доказанным выше утверждением о собственных значениях
и преобразуем \(\sin\theta(x, \lambda_n)\).
\[\sin\theta(x, \lambda_n) = \sin{(\lambda_n^{1/ 2}x+v(x, \lambda_n)+\rho(x, \lambda_n))}\]
Здесь для продолжения цепочки равенств воспользуемся формулой
\eqref{7b}.
\begin{equation*}\begin{aligned}\sin&{(\lambda_n^{1/ 2}x+v(x, \lambda_n)+\rho(x, \lambda_n))} =
\sin{(\mu_n x +mx+v(x, \lambda_n)+\rho(x, \lambda_n))} =
\\ =&\sin(mx)(1+(\mu_n x+v(x, \lambda_n)+\rho(x, \lambda_n))^2) +
\\ &+\cos(mx)(\mu_n x+v(x, \lambda_n)+\rho(x, \lambda_n))
+\rho(x, \lambda_n)= \\  = &\sin(mx) +
\cos(mx)\int_0^x{u(t)\sin{(2mt)}dt}+\frac{1}{2m}\cos(mx)\int_0^x{u^2(t)dt}+\\
&+
2\cos(mx)\int_0^x\int_0^t{u(t)u(s)\cos(2mt)\sin(2ms)dsdt}-\\
&-
\frac{1}{2m}\cos(mx)\int_0^x{u^2(t)\cos(2mt)dt}-\frac{x}{\pi}\cos(mx)\int_0^{\pi}{u(t)\sin{(2mt)}dt}-\\
&-\frac{x}{2m\pi}\cos(mx)\int_0^{\pi}{u^2(t)dt}-\frac{2x}{\pi}\cos(mx)\int_0^{\pi}\int_0^t{u(t)u(s)\cos(2mt)\sin(2ms)dsdt}+\\
& +\frac{x}{2m\pi}\cos(mx)\int_0^{\pi}{u^2(t)\cos(2mt)dt} +
\rho(x, \lambda_n).\end{aligned}\end{equation*} Подставим полученное в
выражения для функций \(y_n(x)\).
\begin{equation*}\begin{aligned}y_n&(x, \lambda_n) = r(x, \lambda_n)\sin\theta(x, \lambda_n) =\sin(mx) + \cos(mx)\int_0^x{u(t)\sin{(2mt)}dt}+\\
&+\frac{1}{2m}\cos(mx)\int_0^x{u^2(t)dt}+2\cos(mx)\int_0^x\int_0^t{u(t)u(s)\cos(2mt)\sin(2ms)dsdt}-
\\
&
-\frac{1}{2m}\cos(mx)\int_0^x{u^2(t)\cos(2mt)dt}-\frac{x}{\pi}\cos(mx)\int_0^{\pi}{u(t)\sin{(2mt)}dt}-\\
&
-\frac{x}{2m\pi}\cos(mx)\int_0^{\pi}{u^2(t)dt}-\frac{2x}{\pi}\cos(mx)\int_0^{\pi}\int_0^t{u(t)u(s)\cos(2mt)\sin(2ms)dsdt}+\\
& +\frac{x}{2m\pi}\cos(mx)\int_0^{\pi}{u^2(t)\cos(2mt)dt}
-\sin(mx)\int_0^x{u(t)\cos(2mt)dt} -\\
& -\frac{1}{2m}\sin(mx)\int_0^x{u^2(t)\sin{(2mt)}dt} + \rho(x,
\lambda_n),\end{aligned}\end{equation*}

при этом \(\{ \sup\limits_{0 \leq x \leq \pi}{|\rho(x,
\lambda_n)|} \} \in l_1\).

Произведем нормировку собственных функций, для этого полученное
выражение необходимо домножить на \(\left(\int_0^{\pi}{y_n(x)
\overline{y_n(x)}dx}\right)^{-1/ 2}\), напомним, что функция
$q(x)$ предполагается комплекснозначной. Будем выделять
вещественную и мнимую часть функции $u(x)$, то есть обзначим $u(x)
= u_R(x) + iu_I(x)$. Тогда сопряженная функция выглядит следующим
образом: $\overline{u(x)} = u_R(x) - iu_I(x)$. Заметим, что при
сложении $y_n(x)$ и $\overline{y_n(x)}$ с одинаковыми
коэффициентами мнимые части будут сокращаться.

При подсчете данного сомножителя важно помнить, что $m$ не
является целым числом, а $m = n - \frac{1}{2}$, где $n$ - целое. В
таком случае, например, \(\int_0^{\pi}{x\sin{(2mx)}} =
\frac{\pi}{2m}\), а не 0.

Итак, после перемножения нормировочный сомножитель имеет следующий
вид:
\begin{multline*}\int_0^{\pi}{y_n(x)\overline{y_n(x)}dx} = \frac{\pi}{2} -
\frac{1}{2m}\int_0^{\pi}{u(t)\sin{(2mt)}dt} - \int_0^{\pi}{(\pi -
t)u_R(t)\cos(2mt)dt}-\\ \ \ - \frac{1}{2m}\int_0^x{(\pi -
t)(u_R^2(t)-u_I^2(t))\sin{(2mt)}dt} + \rho(x, \lambda_n)=
\end{multline*}
\begin{multline*}=
\frac{\pi}{2}\biggl(1 - \frac{1}{\pi m}\int_0^\pi u(t)\sin(2mt)dt
-\frac{2}{\pi}\int_0^\pi (\pi - t)u(t)\cos(2mt)dt-
\\ \ \ -\frac{1}{\pi m}\int_0^x(\pi - t)u^2(t)\sin(2mt)dt\biggr) +
\rho(x, \lambda_n),\end{multline*}
 где $\{ \sup\limits_{0 \leq x \leq \pi}{|\rho(x, \lambda_n)|} \} \in l_1$.

 Рассмотрим подробнее, как получается остаток. Мы не будем приводить подробных выкладок
 возведения в квадрат и взятия интеграла. Разобьем получающиеся слагаемые на типы и
 и отметим, почему каждый из них образует последовательность, лежащую в \(l_1\).
Выделяются следующие типы слагаемых:

 1)
\(\int\limits_0^{\pi}{u(t)\sin{(2mt)}dt}\int\limits_0^{\pi}{u(t)\sin{(2mt)}dt},\)

\noindent изначально \(u(t) \in L_2\), значит ее коэффициенты
Фурье образуют последовательность из пространства \( l_2\), тогда
в квадрате они принадлежат пространству \(l_1\);

2) \(\frac{1}{m}\int\limits_0^{\pi}{u(t)\sin{(2mt)}dt},\)

\noindent аналогично пункту 1) интегральный сомножитель будет
последовательностью из пространства \(l_2\), при домножении его на
\(1/m\), получаем последовательность, принадлежащую \(l_1\);

3)
\(\frac{1}{m}\int\limits_0^{\pi}\int\limits_0^t{u(t)u(s)\cos(2mt)\sin(2ms)dsdt},\)

\noindent последовательность, образуемая двойным интегралом, из
пространства \( l_2\), а сомножитель \(1/m\) переводит
произведение в пространство \(l_1\);

4)\(\frac{1}{m}\int\limits_0^{\pi}{u^2(t)\sin(2mt)dt}\int\limits_0^{\pi}{u(t)\sin(2mt)dt},\)

\noindent второй интеграл образует последовательность,
принадлежащую \(\l_2\), \(u(t) \in L_2\), следовательно, \( u^2(t)
\in L_1  \) и значит ее коэффициенты Фурье убывают; перемножаем
последовательности \(\{\int\limits_0^{\pi}{u(t)\sin(2mt)dt}\}\) и
\(\{1/m\} \) и результирующая последовательность будет из
пространства \(l_1\);

5)
\(\frac{1}{m}\int\limits_0^{\pi}{u^2(t)\sin(2mt)dt}\int\limits_0^{\pi}{u(t)\sin(2mt)dt},\)
\noindent попадает в остаток аналогично пункту 4);

6) все интегралы с сомножителем \(1/m^2\), а так же всевозможные
произведения \(\int\limits_0^{\pi}{u^2(t) dt}\),
\(\int\limits_0^{\pi}{u^2(t)\sin(2mt)dt}\),
\(\int\limits_0^{\pi}{u(t)\sin(2mt) dt}\),
\(\int\limits_0^{\pi}\int\limits_0^t{u(t)u(s)\cos(2mt)\sin(2ms)dsdt}\),
где в формулах встречаются три знака интеграла образуют
последовательности, лежащие в пространстве \(l_1\);

7) все 6 пунктов аналогичны в случаях, когда появляются интегралы
с $\overline{u(x)}$.

Итак осталось возвести полученные выражения для собственных
функций в степень \(-1/ 2\), для чего воспользуемся формулой
Тейлора. В результате получим равенство \eqref{8}. После этого
можем построить и биортогональную систему, воспользовавшись
формулой \(w_n(x) = \frac{\overline{y_n(x)}}{(y_n(x),
\overline{y_n(x)})}\) (см. \cite{S1}), таким образом получим
соотношение \eqref{9b}.

\emph{Теорема доказана.}

\end{document}